\documentclass[a4paper, 11pt, nosumlimits, twoside] {amsart}
\usepackage{amsmath} 
\usepackage{amssymb}
\usepackage{amsfonts}
\usepackage{amsthm}
\usepackage{graphics}
\usepackage{epsfig}
\usepackage{hyperref}
\usepackage[arrow, matrix, curve]{xy}
\makeindex

\numberwithin{subsection}{section}

\newtheoremstyle{note}
  {}
  {}
  {}
  {}
  {\bfseries}
  {.}
  { }
  {}

\newtheoremstyle{notes}
  {}
  {}
  {\itshape}
  {}
  {\bfseries}
  {.}
  { }
  {}

\theoremstyle{notes}
\newtheorem{satz}{Satz}[section]
\newtheorem{lemma}[satz]{Lemma}
\newtheorem{thm}[satz]{Theorem}
\newtheorem{prop}[satz]{Proposition}
\newtheorem{cor}[satz]{Corollary}
\newtheorem{conj}[satz]{Conjecture}

\theoremstyle{note}
\newtheorem{bem}[satz]{Remark}
\newtheorem{fact}[satz]{Fact}
\newtheorem{defn}[satz]{Definition}
\newtheorem{notation}[satz]{Notation}

\newcommand{\R}{\mathbb{R}}														
\newcommand{\N}{\mathbb{N}}														
\newcommand{\Q}{\mathbb{Q}}														
\newcommand{\im}{\textup{im}}													
\newcommand{\Nd}[1]{N^1(#1)_{\mathbb{R}}}							
\newcommand{\Nc}[1]{N_1(#1)_{\mathbb{R}}}							
\newcommand{\num}{\equiv_{\text{num}}}								
\newcommand{\mor}[1]{\overline{\textup{NE}}({#1})}		
\newcommand{\mov}[1]{\overline{\textup{NM}}({#1})}		
\newcommand{\smo}[1]{\overline{\textup{SNM}}({#1})}		
\newcommand{\eff}[1]{\overline{\textup{Eff}}({#1})}		

\newcommand{\ms}{\mapsto}          
\newcommand{\ra}{\rightarrow}      
\newcommand{\Ra}{\Rightarrow}      
\newcommand{\Lra}{\Leftrightarrow} 
\newcommand{\bir}{\dashrightarrow} 

\newcommand{\eps}{\varepsilon}
\newcommand{\ph}{\varphi}
\newcommand{\nb}{\varphi_1^*} 		
\newcommand{\np}{\varphi_{*1}} 		

\begin{document}

\title{The cone of moving curves of a smooth Fano-threefold}

\author{Sammy Barkowski}

\thanks{The author is supported in full by the Graduiertenkolleg
  ``Globale Strukturen in Geometrie und Analysis'' of the Deutsche
  Forschungsgemeinschaft, DFG}
  
\address{Sammy Barkowski, Mathematisches Institut, Universit\"at zu
  K\"oln, Weyertal 86--90, 50931 K\"oln, Germany}
\email{\href{mailto:sbarkows@mi.uni-koeln.de}{sbarkows@mi.uni-koeln.de}}
\urladdr{\href{http://www.mi.uni-koeln.de/~gk/}{http://www.mi.uni-koeln.de/$\sim$gk/}}

\date{April 2, 2007}

\begin{abstract}
In this short note we show that the closed cone of moving curves $\mov{X}\subset\Nc{X}$ of a smooth Fano-threefold $X$ is polyhedral.

The proof relies on a famous result of Bucksom, Demailly, Paun and Peternell which says that the strongly movable cone is dual to the cone of pseudoeffective divisors. Finally, we relate the extremal rays of the cone of moving curves on a threefold with extremal rays in the cone of strongly movable curves on a birational model obtained by a Mori contraction.
\end{abstract}
\maketitle

\setcounter{tocdepth}{1}

\tableofcontents      
\section{Introduction}
In this paper we will give a very simple proof of the fact that the cone $\mov{X}$ of moving curves on a projective variety $X$ is a closed convex polyhedral cone if $X$ is a smooth Fano-threefold:
\begin{thm}\label{mainthm}
Let $X$ be a smooth  Fano-threefold and let $\ph_i:X\bir X_i$, $i=1,...,k$, be the divisorial contractions with exceptional divisors $E_i\subset X$ corresponding to some extremal rays $\R_+[r_i],\ i=1,...,k$, of the Mori cone $\mor{X}$ of $X$. Then $$\mov{X}=\bigl\{[c]\in\mor{X}\mid c\cdot E_i\geq 0\text{, for all }i=1,...,m\bigr\}.$$ In particular, $\mov{X}$ is a closed convex polyhedral cone in $\Nc{X}$.
\end{thm}

This result has been claimed in a much more general setup in the paper \cite{Bat} by Batyrev and the preprint \cite{Arau} by Araujo, but see \cite[Remark~3.4]{Arau}.

To achieve Theorem~\ref{mainthm} we will use a famous result \cite[Theorem~2.2]{bdpp} of Bucksom, Demailly, Paun and Peternell. It says that the strongly movable cone $\smo{X}\subset\Nc{X}$ of a projective variety $X$ is dual to the cone $\eff{X}\subset\Nd{X}$ of pseudoeffective divisors on $X$.

Then we relate the extremal rays of the cone of moving curves on a threefold with extremal rays in the cone of strongly movable curves on a birational model obtained by a Mori contraction.

We have learned that Alex K\"uronya and Endre Szab\'o have independently proved in \cite{Kur} that the cone $\eff{X}$ of pseudoeffective divisors on a projective variety $X$ is a finite rational polyhedron if $X$ is a smooth Fano-threefold. Therefore, one result implies the other by \cite[Theorem~2.2]{bdpp}.
\subsection{Outline of the paper}
In section \ref{cones} we will give the basic definitions and recall the Theorem of BDPP. In section \ref{mmp} we will give a brief review of the minimal model program in dimension 3. We will prove the main claim of this note in section~\ref{moving}. In order to use MMP we want to take ``pullback'' and ``pushforward'' of 1-cycles. Therefore, we have to specify these notions. This will be done in section \ref{pull}.
\subsection{Acknowledgements}\quad The author would like to thank James \hbox{M$^{\textup{c}}$Kernan} for calling his attention to an error in a previous version of the paper. The author is grateful to the DFG and the Graduiertenkolleg ``Glo\-bale Strukturen in Geometrie und Analysis'' for their support. He would also like to thank Stefan Kebekus, Thomas Eckl, Vladimir Lazi\'c and Sebastian Neumann for numerous discussions.
\section{Cones in the N\'eron-Severi spaces}\label{cones}
In the following section let $X$ be a normal projective variety of dimension $n$. Numerical equivalence of $\R$-divisors or 1-cycles will be denoted by $\num$.

\begin{defn}[N\'eron-Severi spaces]
The vector space $\Nc{X}$ of numerical equivalence classes of 1-cycles on $X$ with real coefficients is called the N\'eron-Severi space of curves. The vector space $\Nd{X}$ of numerical equivalence classes of $\R$-divisors on $X$ is called the N\'eron-Severi space of divisors.
\end{defn}
\begin{bem}
The N\'eron-Severi spaces are finite-dimensional real vector spaces. By construction the intersection number of curves and Cartier divisors gives a perfect pairing $$\Nd{X}\times\Nc{X}\ra\R,~(D,c)\ms D\cdot c\in\R.$$ See \cite[Definition~1.4.25]{Laz04}. The dimension $\dim_{\R}(\Nd{X})=:\rho(X)$ of $\Nd{X}$ is called the Picard number of $X$.
\end{bem}\pagebreak
\begin{defn}[Mori Cone]
Let $$\text{NE}(X):=\left\{\sum\limits_{i=1}^m a_i[c_i]\mid c_i\subset X \text{ an irreducible curve}, a_i\geq0\right\}\subset\Nc{X}$$ be the convex cone spanned by all classes of effective 1-cycles on $X$. Its closure $\mor{X}\subset\Nc{X}$ is called the \emph{Mori cone of} $X$.
\end{defn}
\begin{defn}[Movable curve, strongly movable curve]
A class $\gamma\in\Nc X$ is called \emph{movable} or \emph{moving} if there exists an irreducible curve $c\subset X$ with $\gamma=[c]$ such that $c=c_{t_0}$ is a member of an algebraic family $\{c_t\}_{t\in T}\subset X$ of irreducible curves with $\bigcup_{t\in T}c_t=X$.

A class $\gamma\in\Nc{X}$ is called \emph{strongly movable} if there exists a projective birational mapping $\mu:X'\bir X$, together with ample divisors $A_1,...,A_{n-1}$ on $X'$, such that $\gamma=[\mu_*(A_1\cap...\cap A_{n-1})]$.
\end{defn}
\begin{defn}[Moving cone, strongly movable cone]
The closure $$\mov{X}:=\overline{\{\gamma\in\Nc{X}\mid \gamma\text{ a moving class}\}}\subset\Nc{X}$$ of the cone generated by classes of moving curves on $X$ is called the \emph{moving} or \emph{movable cone of} $X$.
The closure $$\smo{X}:=\overline{\{\gamma\in\Nc{X}\mid \gamma\text{ a strongly movable class}\}}\subset\Nc{X}$$ of the cone generated by classes of strongly movable curves on $X$ is called the \emph{strongly movable cone of} $X$.
\end{defn}
\begin{fact}[See \protect{\cite[§2]{bdpp}}]
The cones $\mov{X}$ and $\smo{X}$ are closed convex cones in $\Nc{X}$. We have inclusions $$\smo{X}\subseteq\mov{X}\subseteq\mor{X}.$$
\end{fact}
\begin{defn}[Pseudoeffective cone]
The \emph{pseudoeffective cone} $$\eff{X}\subset\Nd{X}$$ is the closure of the convex cone spanned by the classes of all effective $\R$-divisors on $X$.
\end{defn}
\begin{defn}[Extremal face]\label{face}
Let $K\subset V$ be a closed convex cone in a finite-dimensional real vector space. An \emph{extremal face} $F\subset K$ is a subcone of $K$ having the following property
\begin{equation}\label{extr} \textup{if }v+w\in F\textup{ for some }v,w\in K\textup{, then necessarily }v,w\in F.
\end{equation}
A one-dimensional extremal face $R$ of $K$ is called \emph{extremal ray}.
\end{defn}
\begin{bem}
A subcone with the property \eqref{extr} is often called ``\emph{geometrically extremal}''.
\end{bem}
\begin{defn}[Extremal class]
Let $K$ be closed convex cone in $\Nd{X}$ or in $\Nc{X}$. A class $\gamma\in K$ is called \emph{extremal class} if the ray $r=\R_+\gamma$ spanned by $\gamma$ is an extremal ray of $K$.
\end{defn}
\subsection{Two important theorems of BDPP}
The next theorem is one of the main results in \cite{bdpp}. With this theorem the authors were able to give a characterization of uniruled varieties: a smooth projective variety $X$ is uniruled if and only if $K_X$ is not pseudoeffective.

\begin{thm}[see \protect{\cite[Theorem~2.2]{bdpp}}]\label{bdp}
Let $X$ be an irreducible projective variety of dimension $n$. Then the cones $\smo{X}$ and $\eff{X}$ are dual, i.e. $$\smo{X}=\{\gamma\in\Nc{X}\mid \gamma\cdot D\geq 0,\text{ for all }D\in\eff{X}\}.$$\hfill$\square$
\end{thm}
The following theorem shows why Theorem~\ref{bdp} is useful for our purposes.

\begin{thm}[see \protect{\cite[Theorem~2.4]{bdpp}}]\label{eq}
Let $X$ be an irreducible and smooth projective variety. Then the cones $\smo{X}$ of strongly movable curves and $\mov{X}$  of moving curves coincide.\hfill$\square$
\end{thm}
Therefore, if $X$ is a smooth projective variety, we can use the duality-statement of Theorem \ref{bdp} for our purposes.
\section{The minimal model program in dimension 3}\label{mmp}
In this section we want to give a brief review of the minimal model program in dimension 3.

Even if $X$ is a smooth variety, the minimal model program in dimension 3 can produce singular varieties.
Therefore, we have to introduce the notion of \emph{terminal singularities} as a first step. 

Then we will state the Cone and the Contraction Theorem. Both theorems are essential for the following.
\subsection{The Cone and the Contraction Theorem}
\begin{defn}[Terminal singularities]
A normal variety $X$ of dimension $n$ has only \emph{terminal singularities} if
\begin{enumerate}
\renewcommand{\labelenumi}{\roman{enumi})}
\item $K_X$ is $\Q$-Cartier, i.e. there exists $m\in\N$ such that $mK_X$ is a Cartier divisor,
\item there exists a projective birational morphism $f:Y\bir X$ from a smooth variety $Y$ such 
that all coefficients $a_j$ of the exceptional divisors $E_j$ of $f$ in the ramification formula $K_Y=f^*(K_X)+\sum a_jE_j$ are strictly positive.
\end{enumerate}
\end{defn}
\begin{defn}
A normal variety $X$ is called $\Q$\emph{-factorial} if every Weil divisor $D$ on $X$ is $\Q$-Cartier, i.e. there exists $m\in\N$ such that $mD$ is Cartier.
\end{defn}
Now we are able to state the Cone Theorem.
\begin{thm}[Cone Theorem, see \protect{\cite[Theorem~7-2-1]{Mats}}]\label{conethm}
Let $X$ be a $\Q$-factorial variety with only terminal singularities such that $K_X$ fails to be nef. Then
$$\mor{X}=\mor{X}_{K_X\geq0}+\sum\R_+r_i,$$
where $\mor{X}_{K_X\geq0}:=\{\gamma\in\mor{X}\mid \gamma\cdot K_x\geq0\}$ and the $r_i$ are locally discrete classes in the half-space $\mor{X}_{K_X<0}:=\{\gamma\in\mor{X}\mid \gamma\cdot K_x<0\}$.\\ 
\hspace*{1cm} \hfill$\square$
\end{thm}
\begin{bem}[\protect{\cite[Chapter~6.1~and~7.9]{deb}}]
The proof of the Cone Theorem given in \cite{deb} shows that for every ample $\R$-divisor $A$ and $\eps>0$ there are only finitely many extremal rays in the half space $N_1(X)_{K_X+\eps A<0}$.

This shows that $\mor{X}$ is a convex polyhedral cone if $-K_X$ is ample. See \cite[Example~1.5.34]{Laz04}.
\end{bem}
Another important concept is the notion of an \emph{extremal contraction}. 
\begin{defn}[Extremal contraction]
Let $X$ be a $\Q$-factorial variety with only terminal singularities. Then a morphism $\ph:X\ra Y$ is called an \emph{extremal contraction} with respect to $K_X$ if
\begin{enumerate}
\renewcommand{\labelenumi}{\roman{enumi})}
\item $\ph$ is not an isomorphism,
\item if $c\subset X$ is a curve with $\ph(c)=\textup{pt.}$, then $c\cdot K_X<0$,
\item all the curves contracted by $\ph$ are numerically proportional, i.e. $\ph(c)=\textup{pt.}=\ph(c')\Ra [c]=\lambda [c']$, for a suitable $\lambda\in\Q_+$,
\item $\ph$ has connected fibers with $Y$ being normal and projective.
\end{enumerate}
\end{defn}
\begin{thm}[Contraction Theorem, see \protect{\cite[Theorem~8-1-3]{Mats}}]\label{contr}
Let $X$ be a $\Q$-factorial projective variety with only terminal singularities such that $K_X$ is not nef. For each extremal ray $\R_+r$ of $\mor{X}$ with $r\cdot K_X<0$, there exists a morphism $\ph_r:X\ra Y$, called the \emph{contraction of an extremal ray} $\R_+r$ with respect to $K_X$, such that
\begin{enumerate}
\renewcommand{\labelenumi}{\roman{enumi})}
\item $\ph_r$ is an extremal contraction,
\item if $c\subset X$ is a curve, then $\ph_r(C)=\textup{pt.}\Lra [c]\in\R_+r$.\hfill$\square$
\end{enumerate}
\end{thm}
\begin{defn}
Let $\ph:X\ra Y$ be a birational morphism. Then there exists a birational map $\psi:Y\bir X$ such that $\ph\circ\psi=\textup{id}_Y$ and $\psi\circ\ph=\textup{id}_X$ as rational maps. There exists a maximal open subset $V\subset Y$ such that $\psi|_V:V\ra X$ is a morphism. For details see \cite[Chapter~V.5]{H}. The inverse image $\textup{Exc}_X(\ph):=\ph^{-1}(Y\setminus V)\subset X$ of $Y\setminus V$ under $\ph$ is called the \emph{exceptional locus} of $\ph$ in $X$.
\end{defn}
\begin{bem}
There are three different types of contraction morphisms $\ph:X\ra Y$.
\begin{enumerate}
\item The morphism $\ph$ is birational. If $\textup{codim}_X\textup{Exc}_X(\ph)=1$, we say $\ph$ is a \emph{divisorial contraction} or \emph{of divisorial type}.
\item The morphism $\ph$ is birational. If $\textup{codim}_X\textup{Exc}_X(\ph)\geq2$, we say $\ph$ is a \emph{small contraction} or \emph{of flipping type}.
\item If $\dim Y < \dim X$, we say $\ph$ is \emph{of fiber type}.
\end{enumerate}
In case (1) the exceptional locus $\textup{Exc}_X(\ph)$ consists of a unique irreducible divisor and $Y$ is again $\Q$-factorial with only terminal singularities. The morphism $\ph$ drops the Picard number of $X$ by one, i.e. $\rho(Y)=\rho(X)-1$. In case (3) $X$ is a \emph{Mori fiber space}. In case (2) the morphism $\ph$ is birational, but the canonical divisor $K_Y$ of $Y$ is not $\Q$-Cartier anymore. For details see \cite[Chapters~3~and~8]{Mats}.
\end{bem}

Because of the existence of small contractions, we have to introduce the notion of \emph{flips}.
\begin{conj}
Let $X$ be a $\Q$-factorial projective variety with only terminal singularities and $\ph:X\ra Y$ a small contraction.
\begin{enumerate}
\item There exists a unique birational map $\phi$, from $X$ to a $\Q$-factorial projective variety $X^+$ with only terminal singularities, and another small contraction $\ph^+:X^+\ra Y$ such that the following diagram is commutative
$$\begin{xy}
  \xymatrix{
      X \ar@{-->}[rr]^{\phi} \ar[rd]_{\ph}  &     &  X^+ \ar[dl]^{\ph^+}  \\
                             &  Y  &
  }
\end{xy}
$$
and $K_{X^+}$ is $\ph^+$-ample. We call the map $\phi$ a \emph{flip} of $\ph$.
\item There is no infinite sequence of flips.
\end{enumerate}
\end{conj}

This conjecture is proved in dimension 3, but it is still conjectural in higher dimensions. For details see \cite[Chapter~9]{Mats}.
The minimal model program for higher dimensions stands and falls with this conjecture.
\subsection{Minimal model program}
\begin{enumerate}
\item We start with a $\Q$-factorial projective variety $X$ with only terminal singularities.
\item We ask if $K_X$ is nef. If the answer is yes, we call $X$ a minimal model and stop the program.
If the answer is no, we perform the next step of the program.
\item Since $K_X$ is not nef, there exists a contraction $\ph:X\ra Y$ of an extremal ray $\R_+r$ of $\mor{X}$ with $r\cdot K_X<0$. Now, we ask if $\dim Y<\dim X$. If the answer is yes, $X$ is a Mori fiber space and we stop the program. If the answer is no, we perform the next step of the program.
\item We ask if $\textup{codim}_X\textup{Exc}(\ph)\geq2$. If the answer is no, $Y$ is again $Q$-factorial and projective with only terminal singularities and we restart the program with $Y$ instead of $X$. If the answer is yes, we perform the next step of the program.
\item By the existence of flips, we can flip the map $\ph$ and obtain a $\Q$-factorial projective variety $X^+$ with only terminal singularities. We restart the program with $X^+$ instead of $Y$.
\end{enumerate}

Since the Picard number drops by one with every divisorial contraction and by the termination of flips, this procedure terminates and gives a minimal model or a Mori fiber space as result.
\section{The moving cone of a Fano-threefold}\label{moving}
\begin{defn}
A normal projective variety $X$ of dimension 3 is called a \emph{Fano-threefold} if its anticanonical bundle $-K_X$ is ample.
\end{defn}
The following statement will be the technical basis for our main result.
\begin{prop}\label{prop1}
Let $X$ be a smooth  Fano-threefold and let $\ph_i:X\bir X_i$, $i=1,...,k$, be the divisorial contractions corresponding to some extremal rays $\R_+[r_i],\ i=1,...,k$, of the Mori cone $\mor{X}$ of $X$. If $[r]\in\mor{X}$ is a class with $r\cdot\textup{Exc}_X(\ph_i)\geq 0$, for all $i=1,...,k$, then $r\cdot D\geq 0$ for all irreducible divisors $D\subset X$.
\end{prop}
\begin{proof}
By Theorem \ref{conethm} $\mor{X}$ is a polyhedral convex cone spanned by finitely many extremal rays $\R_+[r_i],\ i=1,...,m$.
Theorem \ref{contr} guarantees the existence of a morphism $$\ph_i:X\ra X_i,$$ 

for every extremal class $[r_i]$ of $\mor{X}$, contracting exactly $\R_+[r_i]$, i.e. contracting all curves $c\subset X$ which are numerically proportional to $r_i$. These contractions are either divisorial or of fiber type, because $X$ is a smooth projective threefold. See \cite[Example~8-3-8]{Mats}. 

Assume that $\ph_i$ is divisorial for $i=1,...,k\leq m$ and let $E_i:=\textup{Exc}_X(\ph_i)\subset X$ denote the exceptional divisor which is contracted by $\ph_i$, for $i=1,...,k$. The divisors $E_i$ are unique and irreducible. See \cite[Proposition~8-2-1]{Mats}.

Now let $[r]\in\mor{X}$ be an arbitrary class with 
\begin{equation}\label{inter}r\cdot E_i\geq 0\text{, for all }i=1,...,k\end{equation} and $[D]\in\eff{X}$ an arbitrary irreducible divisor class on $X$. Because of \eqref{inter}, we can assume that $D\neq E_i$, for all $i=1,...,k$.
Since $\mor{X}$ is polyhedral, we find an effective linear combination of the extremal classes $[r_i]$ such that $$[r]=\sum\limits_{i=1}^m a_i[r_i],\ a_i\geq0.$$ Therefore, to finish the proof it is sufficient to show that $D\cdot r_i\geq 0$, for all $i=1,...,m$.

So let $i\in\{1,...,m\}$.
\begin{enumerate}
\renewcommand{\labelenumi}{Case \arabic{enumi}:}
\item The contraction $\ph_i$, corresponding to the extremal class $[r_i]$, is of fiber type. The divisor $D$ cannot contain all fibers of $\ph_j$. Hence there has to be a fiber $F$ which intersects $D$ properly or $F\cap D=\emptyset$. Since all curves lying in fibers of $\ph_j$ are numerically proportional, we can take a curve $c\in F$ and obtain $0\leq \lambda c\cdot D=r_j\cdot D$, for a suitable $\lambda\in\Q_+$.
\item The contraction $\ph_i$, corresponding to the extremal class $[r_i]$, is divisorial. If $D\cap E_i=\emptyset$, we have $r_i\cdot D=0$, since $r_i$ is contained in $E_i$. If $D\cap E_i\neq\emptyset$, choose a curve $r'\subset E_i$ such that $r'$ is not contained in $D$ and $r'\num \lambda r_i$, for a suitable $\lambda>0$. This is possible, since $D\neq E_i$ and $E_i$ is the unique (!) irreducible divisor containing all curves $c\subset X$ which are numerically proportional to $r_i$. We obtain $0\leq r'\cdot D=\lambda r_i\cdot D$.
\end{enumerate}
This concludes the proof.
\end{proof}
Now we will prove the main result of this note.
\begin{proof}[Proof of Theorem~\protect{\ref{mainthm}}] 
We want to prove that $$\mov{X}=\bigl\{[c]\in\mor{X}\mid c\cdot E_i\geq 0\text{, for all }i=1,...,m\bigr\}.$$
By Theorem \ref{eq} we have $\mov{X}=\smo{X}$, since $X$ is a smooth irreducible projective variety. Theorem \ref{bdp} gives $c\cdot E_i\geq 0$ for all classes $[c]\in\smo{X}=\mov{X}$. Hence $$\mov{X}\subseteq\bigl\{[c]\in\mor{X}\mid c\cdot E_i\geq 0\text{, for all }i=1,...,m\bigr\}.$$
Proposition~\ref{prop1} and Theorem~\ref{bdp} give $$\mov{X}\supseteq\bigl\{[c]\in\mor{X}\mid c\cdot E_i\geq 0\text{, for all }i=1,...,m\bigr\}.$$ 
\end{proof}
\begin{figure}[htbp]
  \centering
  \includegraphics[width=6.3cm]{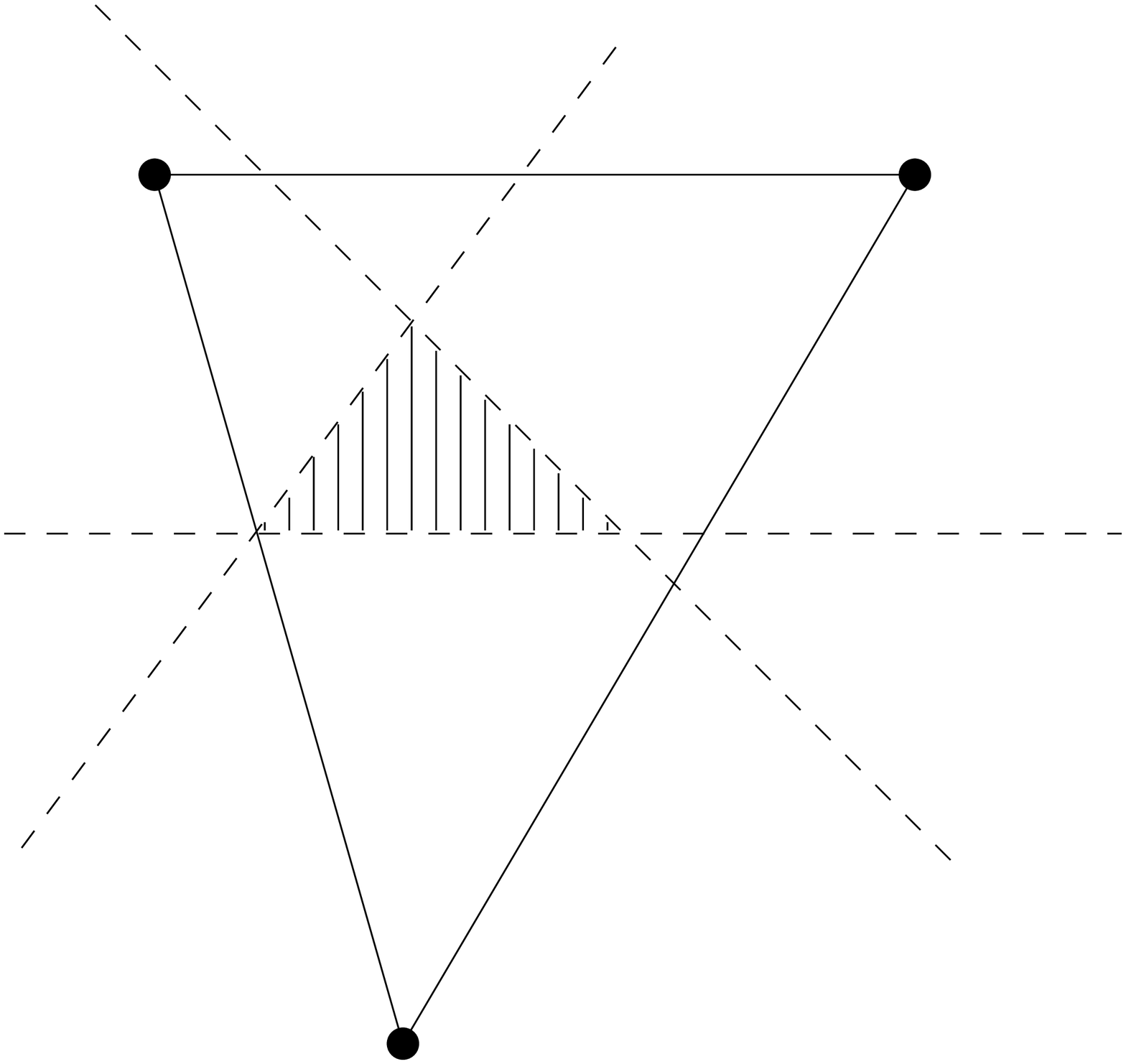}%
		\begin{picture}(0,0)%
				\put(-34.0,-.5){cross-section through $\Nc{X}$}
				\put(-97.0,103.0){$\mov{X}$}
				\put(-93.0,28.0){$\mor{X}$}
				\put(-47.0,147.0){$r_1$}
				\put(-133.0,-.5){$r_2$}
				\put(-170.0,147.0){$r_3$}
				\put(-41.0,50.0){$E_1^{\bot}$}
				\put(-21.0,90.0){$E_2^{\bot}$}
				\put(-101.0,160.0){$E_3^{\bot}$}
		\end{picture}%
  \caption{The hyperplanes~$E_i^{\bot}=\{\gamma\in\Nc{X}\mid \gamma\cdot E_i=0\}$, $i=1,2,3$, cut out the moving cone $\mov{X}$ of $X$.}
  \label{fig:cone} 
\end{figure}
\begin{cor}
In the situation of Proposition \ref{prop1} the following holds. Let $[c]$ be an arbitrary extremal class of $\mov{X}$. Then there exists a divisorial contraction $\ph:X\bir X'$, corresponding to an extremal class $r$ of $\mor{X}$, such that $c\cdot\textup{Exc}_X(\ph)=0$.\hfill $\square$
\end{cor}
\section{Numerical pullback of extremal classes}\label{pull}
\subsection{Basic definitions and notation}
This section is a short survey of the notions given in \cite[Section~3]{Arau}. A detailed treatment of pushforward and pullback of $k$-cycles is given in \cite{Ful}. 
\begin{notation}
Let $X$ and $Y$ be $\Q$-factorial varieties with only terminal singularities and $\ph:X\ra Y$ be a proper morphism. 
\end{notation}
\begin{fact}\label{fact}
The pullback of $\Q$-Cartier divisors on $Y$ gives an injective linear map $$\ph^*:\Nd{Y}\hookrightarrow\Nd{X}$$ and the pushforward of $\Q$-Cartier divisors on $X$ gives a surjective map $$\ph_*:\Nd{X}\twoheadrightarrow\Nd{Y}$$ such that $\ph_*\circ\ph^*=\textup{id}_{\Nd{Y}}$. See \cite[Definition~3.1]{Arau}.
\end{fact}
\begin{defn}[Numerical pullback and pushforward]
Let $$\nb:\Nc{Y}\hookrightarrow\Nc{X}$$ be the dual linear map of the pushforward $\ph_*:\Nd{X}\twoheadrightarrow\Nd{Y}$ of $\R$-divisors on $X$ and let $$\np:\Nc{X}\twoheadrightarrow\Nc{Y}$$ be the dual linear map of the pullback $\ph^*:\Nd{Y}\hookrightarrow\Nd{X}$ of $\R$-divisors on $Y$. We call $\nb$ the \emph{numerical pullback via} $\ph$ and $\np$ the \emph{numerical pushforward via} $\ph$.
\end{defn}
\begin{bem}\label{linmap}
Every 1-cycle $c\in\Nc{Y}$ on $Y$ can be considered as a linear map $$c:\Nd{Y}\ra\R,\ D\ms c\cdot D.$$ 

Therefore, the 1-cycle $\nb(c)\in\Nc{X}$ can be considered as the composition $$c\circ\ph_*:\Nd{X}\ra\R,\ D'\ms c\cdot\ph_*(D).$$

In the same way $\np(c')\in\Nc{Y}$ can be considered as the map $$c'\circ\ph^*:\Nd{Y}\ra\R,\ \tilde{D}\ms c'\cdot\ph^*(\tilde{D})$$ 

for every $c'\in\Nc{X}$. This yields the following two statements.
\end{bem}
\begin{lemma}[Projection formula]\label{proj}
Let $\ph:X\ra Y$ be a proper morphism of complete varieties or projective schemes.
\begin{enumerate}
\renewcommand{\labelenumi}{\roman{enumi})}
\item If $c\in\Nc{Y}$ and $D\in\Nd{X}$, then $\nb(c)\cdot D=c\cdot\ph_*(D)$.\label{pro2}
\item If $c\in\Nc{X}$ and $D\in\Nd{Y}$, then $c\cdot\ph^*(D)=\np(c)\cdot D$.
\end{enumerate}\hfill$\square$
\end{lemma}
\begin{lemma}[\protect{\cite[Remark~3.2]{Arau}}]\label{rmk}
The composition $\np\circ\nb$ is the identity map on $\Nc{Y}$ and
the numerical pullback via $\ph$ satisfies the following conditions.
\begin{enumerate}
\renewcommand{\labelenumi}{\roman{enumi})}
\item If $D\in\Nd{Y}$ and $c\in\Nc{Y}$, then $\ph^*(D)\cdot\nb(c)=D\cdot c$.
\item If $D\in\ker\ph_*$ and $c\in\im\nb$, then $D\cdot c=0$.
\end{enumerate}\hfill$\square$
\end{lemma}
\subsection{Numerical pullback of strongly movable extremal classes}
Finally, we want to propose a method to compute extremal rays of the moving cone of a smooth projective variety. It is based on ideas of Carolina Araujo. In \cite{Arau} she takes numerical pullback of curves lying in general fibers of Mori fiber spaces obtained by running the minimal model program.

We want to run the minimal model program and take the numerical pullback of extremal rays of the moving cones of the obtained varieties successively. Since the varieties obtained can be singular, we have to consider the cone of strongly movable curves instead of the moving cone. Therefore, we show a correspondence between the extremal rays of the strongly movable cone of a $\Q$-factorial projective variety and extremal rays of the strongly movable cone of a birational model obtained by a divisorial contraction.
\begin{thm}
Let $X$ be a $\Q$-factorial projective variety with only terminal singularities and $\ph:X\bir Y$ a divisorial contraction of an extremal ray of $\mor{X}$. Then the following holds.

A class $c\in\Nc{Y}$ is an extremal class of $\smo{Y}$ if and only if the class $\nb(c)\in\Nc{X}$ is an extremal class of $\smo{X}$.
\end{thm}
\begin{proof}[Proof. Step 1.]
Let $c\in\Nc{Y}$ be an extremal class of $\smo{Y}$. We want to prove that $\nb(c)\in\Nc{X}$ is an extremal class of $\smo{X}$.

Let $D\in\eff{X}$ be a effective divisor class on $X$. The projection formula \ref{proj} gives $$\nb(c)\cdot D=c\cdot \ph_*(D).$$ 

Since $\ph_*(D)$ is again effective and $c\in\smo{Y}$, Theorem \ref{bdp} says $\nb(c)\cdot D\geq0$ and hence $\nb(c)\in\smo{X}$.

Now let $v,w\in\smo{X}$ be arbitrary classes with $v+w\in\R_+\nb(c)$ and let $E:=\textup{Exc}_X(\ph)\subset X$ be the exceptional divisor, which is contracted by $\ph$. Thanks to Lemma~\ref{rmk}~\textup{ii)}, we have $$0=\nb(c)\cdot E=\lambda(v+w)\cdot E,$$ for a suitable $\lambda\in\R_+$. Therefore, $v\cdot E = -w\cdot E$ and since $E$ is effective, again by Theorem~\ref{bdp} we find $0\leq v\cdot E=-w\cdot E\leq 0\Ra v\cdot E=0=w\cdot E$. This implies 
\begin{equation}\label{id}
v=\nb(\np(v))\text{ and }w=\nb(\np(w)).
\end{equation}

Now let $D'\in\eff{Y}$ be an effective divisor class on $Y$. Then $\ph^*(D')$ is again effective. Together, the projection formula and Theorem \ref{bdp} give $$\np(v)\cdot D'=v\cdot\ph^*(D')\geq0\text{ and }\np(w)\cdot D'=w\cdot\ph^*(D')\geq0.$$

Hence $\np(v),\np(w)\in\smo{Y}$.

If $\hat{D}\in\Nd{Y}$ is a divisor class on $Y$, we find 
\begin{align*}
(\np(v)+\np(w))\cdot\hat{D}&=\np(v+w)\cdot\hat{D}\\
&\stackrel{\ref{proj}}{=}(v+w)\cdot\ph^*(\hat{D})\\
&=\lambda\nb(c)\cdot\ph^*(\hat{D})\\
&\stackrel{\ref{rmk}\textup{ i)}}{=}\lambda c\cdot \hat{D}.
\end{align*}
Hence $(\np(v)+\np(w))\in\R_+c$. Since $c\in\smo{Y}$ is an extremal class, $\np(v),\np(w)\in\R_+c$. Now, for suitable $\mu,\mu'\in\R_+$ and $\tilde{D}\in\Nd{X}$ an arbitrary divisor class on $X$, the following holds:
$$\nb(c)\cdot\tilde{D}\stackrel{\ref{proj}}{=}c\cdot\ph_*(\tilde{D})=\mu\np(v)\cdot\ph_*(\tilde{D})
\stackrel{\ref{proj}}{=}\mu\nb(\np(v))\cdot\tilde{D}\stackrel{\eqref{id}}{=}\mu v\cdot\tilde{D},$$
$$\nb(c)\cdot\tilde{D}\stackrel{\ref{proj}}{=}c\cdot\ph_*(\tilde{D})=\mu'\np(w)\cdot\ph_*(\tilde{D})
\stackrel{\ref{proj}}{=}\mu'\nb(\np(w))\cdot\tilde{D}\stackrel{\eqref{id}}{=}\mu' w\cdot\tilde{D}.$$

This implies $v,w\in\R_+\nb(c)$. So, $\nb(c)\in\smo{X}$ is an extremal class.

\emph{Step 2 of the proof.} Now let $\nb(c)\in\smo{X}$ be an extremal class for some class $c\in\Nc{Y}$. We want to prove that $c\in\Nc{Y}$ is an extremal class of $\smo{Y}$.

Let $D\in\eff{Y}$ be an effective divisor class on $Y$. Together, Lemma~\ref{rmk}~\textup{ii)} and Theorem \ref{bdp} give 
$$c\cdot D=\nb(c)\cdot\ph^*(D)\geq0,$$ since $\ph^*(D)\in\Nd{X}$ is again effective. Theorem~\ref{bdp} again implies $c\in\smo{Y}$.

We choose arbitrary classes $v,w\in\smo{Y}$ with $v+w\in\R_+c$. If $D'\in\eff{X}$ is an effective divisor class on $X$, then $\ph_*(D')\in\Nd{Y}$ is effective, too. As before, the projection formula and Theorem \ref{bdp} guarantee
$$\nb(v)\cdot D'=v\cdot\ph_*(D')\geq0\text{ and }\nb(w)\cdot D'=w\cdot\ph_*(D')\geq0.$$

This implies $\nb(v),\nb(w)\in\smo{X}$.

For a suitable $\lambda\in\R_+$ and an arbitrary divisor class $\hat{D}\in\Nd{X}$, the following holds:
\begin{align*}
(\nb(v)+\nb(w))\cdot\hat{D}&=\nb(v+w)\cdot\hat{D}\\
&\stackrel{\ref{proj}}{=}(v+w)\cdot\ph_*(\hat{D})\\
&=\lambda c\cdot\ph_*(\hat{D})\\
&\stackrel{\ref{proj}}{=}\lambda\nb(c)\cdot\hat{D}.
\end{align*}
We find $(\nb(v)+\nb(w))\in\R_+\nb(c)$ and $\nb(c)$ is an extremal class of $\smo{X}$ by assumption. This implies $\nb(v),\nb(w)\in\R_+\nb(c)$ and enables us to conclude the proof.

Let $\tilde{D}\in\Nd{Y}$ be an arbitrary divisor class. For suitable $\mu,\mu'\in\R_+$, the following holds:
$$v\cdot\tilde{D}\stackrel{\ref{rmk}\textup{ i)}}{=}\nb(v)\cdot\ph^*(\tilde{D})
=\mu\nb(c)\cdot\ph^*(\tilde{D})\stackrel{\ref{rmk}\textup{ i)}}{=}\mu c\cdot\tilde{D},$$
$$w\cdot\tilde{D}\stackrel{\ref{rmk}\textup{ i)}}{=}\nb(w)\cdot\ph^*(\tilde{D})
=\mu'\nb(c)\cdot\ph^*(\tilde{D})\stackrel{\ref{rmk}\textup{ i)}}{=}\mu' c\cdot\tilde{D}.$$

Hence $v,w\in\R_+c$ and $c\in\smo{Y}$ is an extremal class.
\end{proof}
\def\cprime{$'$}


\begin{thebibliography}{BDPP04}

\bibitem[Ara05]{Arau}
Carolina Araujo.
\newblock {\em The cone of effective divisors of log varieties after
  {B}atyrev}.
\newblock math.AG/0502174, 2005.
\newblock Preprint.

\bibitem[Bat92]{Bat}
V.V. Batyrev.
\newblock {\em The cone of effective divisors of threefolds}.
\newblock Proceedings of the {I}nternational {C}onference on {A}lgebra, {Part}
  3 ({Novosibirsk, 1989}). Contemp. Math., vol. 131, AMS, 1992.
\newblock pp. 337-352.

\bibitem[BDPP04]{bdpp}
Sebastien Boucksom, Jean-Pierre Demailly, Mihai Paun, and Thomas Peternell.
\newblock {\em The pseudo-effective cone of a compact {K}\"ahler manifold and
  varieties of negative {K}odaira dimension}.
\newblock 2004.

\bibitem[Deb01]{deb}
Olivier Debarre.
\newblock {\em Higher-dimensional algebraic geometry}.
\newblock Universitext. Springer-Verlag, New York, 2001.

\bibitem[Ful98]{Ful}
William Fulton.
\newblock {\em Intersection theory}, volume~2 of {\em Ergebnisse der Mathematik
  und ihrer Grenzgebiete. 3. Folge. A Series of Modern Surveys in Mathematics
  [Results in Mathematics and Related Areas. 3rd Series. A Series of Modern
  Surveys in Mathematics]}.
\newblock Springer-Verlag, Berlin, second edition, 1998.

\bibitem[Har77]{H}
Robin Hartshorne.
\newblock {\em Algebraic geometry}.
\newblock Springer-Verlag, New York, 1977.
\newblock Graduate Texts in Mathematics, No. 52.

\bibitem[KS]{Kur}
Alex K\"uronya and Endre Szab\'o.
\newblock {\em The locally polyhedral structure of the effective cone}.
\newblock Unpublished Preprint.

\bibitem[Laz04]{Laz04}
Robert Lazarsfeld.
\newblock {\em Positivity in algebraic geometry. {I}}, volume~48 of {\em
  Ergebnisse der Mathematik und ihrer Grenzgebiete. 3. Folge. A Series of
  Modern Surveys in Mathematics [Results in Mathematics and Related Areas. 3rd
  Series. A Series of Modern Surveys in Mathematics]}.
\newblock Springer-Verlag, Berlin, 2004.
\newblock Classical setting: line bundles and linear series.

\bibitem[Mat02]{Mats}
Kenji Matsuki.
\newblock {\em Introduction to the {M}ori program}.
\newblock Universitext. Springer-Verlag, New York, 2002.

\end{thebibliography}
\end{document}